\documentclass[12pt,leqno]{article}
\usepackage{amsmath, amsthm, amsfonts, amssymb, color}
\usepackage{mathrsfs}
\theoremstyle{definition}

\newcommand{\scr}[1]{\mathscr #1}
\definecolor{wco}{rgb}{0.5,0.2,0.3}

\numberwithin{equation}{section} \theoremstyle{remark}

\title{{\bf  Poincar\'e Inequality on the Path Space of Poisson Point Processes}\footnote{Supported in part by WIMCS,
Creative Research Group Fund of
 the National Natural Science Foundation of
China (No. 10721091) and the 973-Project.} }
\author{{\bf  Feng-Yu Wang$^{a,b}$\thanks{ Corresponding author. E-mail: wangfy@bnu.edu.cn, F.Y.Wang@swansea.ac.uk, C.Yuan@swansea.ac.uk }\,  and\, Chenggui Yuan$^{b}$}\\
\footnotesize{$^a$ School  of Mathematical Sciences \&  Lab. Math. Com. Sys.,} \\\footnotesize{Beijing Normal
University, Beijing 100875, China}\\
 \footnotesize{$^b$  Department of Mathematics, Swansea University, Singleton Park, SA2 8PP, Swansea, UK}}
\begin{document}
\maketitle
\def\R{\mathbb R} \def\Z{\mathbb Z} \def\ff{\frac} \def\ss{\sqrt}
\def\N{\mathbb N}
\def\dd{\delta} \def\DD{\Delta} \def\vv{\varepsilon} \def\rr{\rho}
\def\<{\langle} \def\>{\rangle} \def\GG{\Gamma} \def\gg{\gamma}
\def\ll{\lambda} \def\LL{\Lambda} \def\nn{\nabla} \def\pp{\partial}
\def\d{{\rm d}} \def\bb{\beta} \def\aa{\alpha} \def\D{\scr D}
\def\E{\scr E} \def\si{\sigma} \def\ess{\text{\rm{ess}}}
\def\beg{\begin} \def\beq{\begin{equation}}  \def\F{\scr FC}
\ \def\OO{\Omega}\def\e{{\rm e}}\def\o{\theta}
\def\oo{\omega}  \def\cut{\text{\rm{cut}}} \def\P{\mathbb P} \def\Z{\mathbb Z}
\def\g{\gg} \def\mb{\mathbb} \def\Q{\mathbb Q}\def\B{\scr B}\def\tt{\tilde}
\date{}

\maketitle

\begin{abstract} The quasi-invariance is proved for the distributions of Poisson
point processes under a random shift map on the path space. This leads to a
natural Dirichlet form of jump type on the path space. Differently
from the O-U Dirichlet form on the Wiener space satisfying the
log-Sobolev inequality, this Dirichlet form  merely satisfies the
Poincar\'e inequality but not the log-Sobolev one.

\end{abstract}

\section {Introduction}

Stochastic analysis on the path space of diffusion processes on
manifolds has been intensively studied in the past 15 years, see
e.g. \cite{D} for the quasi-invariance of (Riemannian) Wiener
measures and integration by parts for the Malliavin gradient,
\cite{F, AE, Hsu, CHL} for Poincar\'e and log-Sobolev inequalities
of the associated O-U type Dirichlet forms, and \cite{Feyel, W04,
FW, DGW} for Talagrand type transportation cost inequalities,
\cite{DR, L} for construction of the associated
infinite-dimensional diffusion processes, and \cite{EL1, EL2} for
the study of $L^2$-Hodge theory and the Markovian uniqueness of the
Dirichlet form.

The purpose of this note is to start the corresponding analysis on
the path space of jump processes. As a standard model, we shall
consider here the Poisson process $X=\{X_t:\ 0\le t\le T\}$ on
$\R^d$  with $X_0=0$ and the intensity $\nu$ (cf. \cite{O05}). We assume that
$\nu$ is a   probability
measure on $\R^d\setminus\{0\}$ . The study for finite  $\nu$ is
equivalent to a change of the time interval $[0,T]$.

   The path space of $X$ is

$$M_T=\{\oo: [0,T]\to \R^d:\ \oo\ \text{is\ right\ continuous\ having \
left\ limits}\},$$ which is a Polish space under the Skrohod metric

\beg{equation*}\beg{split} \rr(\oo,\gg):= \inf &\Big\{\dd>0:\
\text{there\ exist\ }n\ge 1, 0=s_0<\cdots<s_n=T, 0=t_0<\cdots
<t_n=T\\
& \text{such\ that}\ |t_i-s_i|\lor |\oo_T-\gg_T|\le \dd, \sup_{s\in
[s_{i-1}, s_i), t\in [t_{i-1}, t_i)}|\oo_s-\gg_t|\le \dd,\ 1\le i\le
n\Big\}.\end{split}\end{equation*} Let $\mu$ be the distribution of
$X$, which is thus a probability measure on $M_T$.

Following the definition of the Malliavin gradient on the Wiener
space, we need to introduce a shift operator on $M_T$ such that
$\mu$ is quasi-invariant. Intuitively, if the process jumps at time
$t$, the path from $t$ on should be shifted in the same scale. On
the other hand, however, since in probability one the process does
not jump at a fixed time, to make $\mu$ quasi-invariant the path can
only be shifted at a random time. As a simple choice, we shall take
below the random time an independent uniform random variable $\tau$
on $[0,T].$ Thus, the map will be

$$\psi: \oo \mapsto \oo + \xi 1_{[\tau, T]},$$ where $\xi$ is a
random variable on $\R^d$ with distribution $\nu$ such that $\{X,
\tau, \xi\}$ are independent.

Let $\tt\mu$ be the distribution of $X+\xi 1_{[\tau, T]}$, which is
again a probability measure on $M_T$. According to Theorem \ref{T2.2} below we have

\beq\label{1.1} \tt\mu(\d\oo)= T^{-1}
N_T(\oo)\mu(\d\oo),\end{equation} where under the probability
measure $\mu$,

$$N_t(\oo):=\sum_{s\le t} 1_{\{\oo_s\ne \oo_{s-}\}},\ \ \ t\ge 0$$
is a Poisson process on $\Z_+$ with intensity $\dd_1$ (i.e. the Dirac measure at $1$).  In particular, 
$\mu(N_T)=T.$

For two functions $F,G$ on $M_T$, let

$$\GG_{t, x} (F,G)= \big(F(\cdot+
x1_{[t,T]})-F(\cdot)\big)\big(G(\cdot+x1_{[t,T]})-G(\cdot)\big),\ \ \ x\in
\R^d, t\in [0,T].$$ Due to (\ref{1.1}),  the form

$$\E(F,G):= \ff 1 T\int_{M_T} \mu(\d\oo) \int_0^T\d t\int_{\R^d}
\GG_{t,x}(F,G)(\oo)\nu(\d x)$$ is well defined on

$$\D(\E):= \big\{F\in L^2(M_T; \mu):\ \E(F,F)<\infty\big\};$$ that
is, the value of $\E(F,G)$ does not depend on $\mu$-versions of $F$
and $G$. The main result of the paper is the following.

\beg{thm}\label{T1.1} $(\E,\D(\E))$ is a symmetric, conservative
Dirichlet from on $L^2(M_T;\mu)$ and the Poincar\'e inequality

\beq\label{P} \mu(F^2)\le T\E(F,F)+ \mu(F)^2,\ \ \ F\in
\D(\E)\end{equation}  holds. But for any $C>0$ the log-Sobolev
inequality

\beq\label{LS} \mu(F^2\log F^2)\le C \E(F,F),\ \ \ F\in \D(\E),
\mu(F^2)=1\end{equation} does not hold. \end{thm}

This result will be proved in the next two sections:  (\ref{1.1}) will be proved in Section 2 while the remainder, i.e. the proof of (\ref{P}) and the disproof of (\ref{LS}), will be 
 addressed in Section 3.

\section{The Dirichlet form and generator}

The main purpose of this section is to prove (\ref{1.1}). We shall first prove 
it for Markov chains on $\Z^d$ then extend to Poisson processes by an approximation argument. 

For any $k\in \Z^d$, let $$N_T^{(k)}(\oo)= \sum_{t\le T} 1_{\{\oo_t- \oo_{t-}=k\}},\ \ \oo \in M_T.$$ 

\beg{lem}\label{L2.1} Let $\nu$ be supported on $\Z^d$.  Then for any $k\in \Z^d$ such that $\nu(k)>0$, the distribution $\mu_k$ of $X+k1_{[\tau, T]}$ is absolutely continuous with respect to $\mu$ with

\beq\label{2.0} \ff{\d\mu_k}{\d \mu}(\oo)= \ff 1 {T \nu(k)} N_T^{(k)}(\oo).\end{equation} Consequently, 
$(\ref{1.1})$ holds and  $(\E,\D(\E))$ is a well-defined 
Dirichlet form on $L^2(M_T;\mu).$\end{lem}

\beg{proof} We shall only prove (\ref{2.0}) since the proof of the remainder is simple and standard. In the present framework $X$ is a shift-invariant conservative  Markov chain on $\Z^d$ with the $Q$-matrix $(q_{ij})_{i,j\in \Z^d}$ given by (cf. \cite{Chen})

$$q_{ij} = \beg{cases} \nu(i-j), &\text{if}\ i\ne j, \\
-\sum_{k\ne i} \nu(k-i), &\text{if}\ i=j.\end{cases}$$ Let $\{p_t(l):\ l\in \Z^d\}$ be the distribution of $X_t$. On any finite set $K$ there exist a constant $C>0$ and a  positive continuous function $h$ on $[0,1]$ with $h(s)\downarrow 0$ as $s\downarrow 0$ such that 

\beq\label{2.2} |p_s(0)-1|\le C s,\ \ \ |p_s(l)-s\nu(l)|\le h(s)s,\ \ l\in K\setminus\{0\}, s\in [0,1].\end{equation} Let $F$ be a bounded  cylindrical function on $M_T$ depending only on coordinates at $0<s_1<\cdots<s_m\le T.$ It suffices to show that 

\beq\label{2.31} \mathbb EF(X+k1_{[\tau,T]})= \ff 1 {T\nu(k)} \mathbb E [F(X) N_T^{(k)}].\end{equation} By first replacing $\nu(k)$ by $\nu_\vv(k):= (\nu(k)+ \vv 2^{-|k|})/(1+\vv \sum_k 2^{-|k|})$ then letting $\vv\downarrow 0$, we may and do assume that $\nu(k)>0$ for all $k\in \Z^d.$ 

To prove (\ref{2.31}), let $I_n=\{0=t_0<t_1<\cdots<t_n=T\}\supset \{s_1,\cdots, s_m\}$ be a sequence of partitions such that $\dd(I_n):= \max_{1\le i\le n} (t_i-t_{i-1})\downarrow 0$. Let $\|\oo\|_n= \sup\{|\oo_{t_i}|:\ t_i\in I_n\}.$ Let $f$ be a bounded measurable function on $(\Z^d)^n$ such that 

$$F(\oo)= f(\oo_{t_1},\cdots, \oo_{t_n}),\ \ \ \oo\in M_T.$$ Let $l_0=0$. For any $R>0$, we have

\beg{equation*}\beg{split} U_R&:= \mathbb E\big[ 1_{\{\|X\|_n\le R\}}F(X+k1_{[\tau,T]})\big]=\ff 1 T \int_0^T \mathbb E\big[  1_{\{\|X\|_n\le R\}}F(X+k1_{[t,T]})\big]\d t\\
&= \ff 1 T \sum_{i=1}^n (t_i-t_{i-1}) \mathbb E\big[  1_{\{\|X\|_n\le R\}}F(X+k1_{[t_i,T]})\big]\\
&=
 \ff 1 T \sum_{i=1}^n (t_i-t_{i-1}) \mathbb E\big[  1_{\{\|X\|_n\le R\}}f(X_{t_1},\cdots, X_{t_{i-1}}, X_{t_i}+k,\cdots, X_{t_n}+k)\big]\\
 &= \ff 1 T \sum_{|l_1|,\cdots, |l_n|\le R}\sum_{i=1}^n (t_i-t_{i-1}) f(l_1,\cdots, l_{i-1}, l_i +k,\cdots, l_n +k) \prod_{j=1}^np_{t_j-t_{j-1}}(l_j-l_{j-1}).
 \end{split}\end{equation*}
By changing variables $(l_1,\cdots, l_n)\mapsto (l_1,\cdots, l_{i-1}, l_i -k, \cdots, l_n-k)$,  and letting

$$K_{R,k,i}=\big\{(l_1, \cdots, L_n)\in (\Z^{d})^n:  |l_j|\le R\text{\ for\ }j\le i-1, |l_j-k|\le R\ \text{for}
\ j\ge i\big\},$$
we obtain 

\beg{equation*}\beg{split} U_R&= \ff 1 T \sum_{i=1}^n (t_i-t_{i-1}) \sum_{(l_1,\cdots l_n)\in K_{R,k,i}} f(l_1,l_2,\cdots, l_n)\ff{p_{t_i-t_{i-1}}(l_i-l_{i-1}-k)}{p_{t_i-t_{i-1}}(l_i-l_{i-1})}
\prod_{j=1}^np_{t_j-t_{j-1}}(l_j-l_{j-1})\\
&= \ff 1 T \sum_{i=1}^n(t_i-t_{i-1})\mathbb E\Big[ 1_{K_{R,k,i}}(X_{t_1},\cdots, X_{t_n}) F(X) \ff{p_{t_i-t_{i-1}}(X_{t_i}-X_{t_{i-1}}-k)}{p_{t_i-t_{i-1}}(X_{t_i}-X_{t_{i-1}})}\Big].\end{split}\end{equation*} By this and (\ref{2.2}), for small enough $\dd(I_n)$ such that 
(note that we have assumed that $\nu(l)>0$ for all $l\in \Z^d$)
$$
\min\{\nu(l): |l|\le k+2R\}-C h \circ\dd(I_n)>0,
$$
where $C$ and $h$ are defined in \eqref{2.2} for $K= \{l\in \Z^d:\ |l|\le k+2R\},$   we have 

\beg{equation*}\beg{split} &\Big|U_R- \ff 1 T \sum_{i=1}^n (t_i-t_{i-1}) 
\mathbb E\Big[1_{K_{R,k,i}}(X_{t_1},\cdots, X_{t_n})F(X) 1_{\{X_{t_i}-X_{t_{i-1}}=k\}}
 \ff{p_{t_i-t_{i-1}}(0)}{p_{t_i-t_{i-1}}(k)}\Big]\Big|\\
 &\le \ff{\|F\|_\infty} T \sum_{i=1}^n (t_i-t_{i-1}) \ff{\max\{\nu(l):\ |l|\le |k|+2R\} + h\circ\dd(I_n)}{\min\{\nu(l): |l| \le k+2R\}- h\circ \dd(I_n)}\mathbb P\big(
 X_{t_i}-X_{t_{i-1}}\notin\{k,0\} \big)\\
 &\quad +\ff{\|F\|_\infty\dd(I_n)[\nu(k)+C h\circ\dd(I_n)]} {T(1-C\dd(I_n))} \sum_{i=1}^n (t_i-t_{i-1}) \mathbb P(X_{t_i}=X_{t_{i-1}})\\
 &\le  \ff{C'\|F\|_\infty \dd(I_n)} T\mathbb E N_T+ \ff{\|F\|_\infty \dd(I_n)[(\nu(k) + C h\circ\dd(I_n)]} {1-C\dd(I_n)}\end{split}\end{equation*}
 for some  $C'>0$  depending only on $\nu$ and $|k|+2R.$ Here, we have used the fact that
 
 $$\sum_{i=1}^n \mathbb P(X_{t_i}\ne X_{t_{i-1}})\le \mathbb E\sum_{i=1}^n 1_{\{X_{t_i}\ne X_{t_{i-1}}\}}
 \le \mathbb E N_T.$$ Letting $n\to\infty$ and using (\ref{2.2}) again, we arrive at 
 
 \beg{equation*}\beg{split} U_R &= \lim_{n\to\infty} \ff 1 T \sum_{i=1}^n (t_i-t_{i-1}) \mathbb
 E\Big[1_{K_{R,k,i}}(X_{t_1},\cdots, X_{t_n})F(X) 1_{\{X_{t_i}-X_{t_{i-1}}=k\}}
 \ff{p_{t_i-t_{i-1}}(0)}{p_{t_i-t_{i-1}}(k)}\Big]\\
 &=\ff 1 {T\nu(k)} \mathbb E\big[F(X)N_T^{(k)}1_{\{\|X-k1_{[\tau,T]}\|_\infty\le R\}}\big],
 \end{split}\end{equation*}where $\|\oo\|_\infty:= \sup_{t\le T} |\oo_t|.$ Letting $R\to\infty$ we complete the proof.\end{proof}

To identify the generator of $(\E,\D(\E))$,  for any $k\in \Z^d$ let $\pi_k(\cdot, \d t)$ be the regular conditional distribution of $\tau$ given $X+k1_{[\tau,T]}.$

\beg{lem}\label{L2.2}  If $\nu$ is supported on $\Z^d$ then the generator $(L,\D(L))$ of $(\E,\D(\E))$ satisfies $\D(L)\supset \B_b(M_T)$, the class of all 
bounded measurable functions on $M_T$, and 

\beq\label{2.1} \beg{split} L F(\oo)= &\sum_{k\in \Z^d} \ff{\nu(k)}T \int_0^T \Big(F(\oo+ k1_{[t,T]})
-F(\oo)\Big)\d t\\
&+\sum_{k\in \Z^d} \ff{N_T^{(k)}(\oo)}T\int_0^T\Big(F(\oo- k1_{[t,T]})-F(\oo)\Big)\pi_k(\oo,\d t),\ \ \ F\in \B_b(M_T).\end{split}\end{equation}\end{lem} 

\beg{proof} Let $F,G\in \B_b(M_T).$  We have 

\beg{equation}\label{2.32}\beg{split} -\E(F,G) &=\sum_{k\in \Z^d} \ff{\nu(k)}T \int_0^T \d t \int_{M_T} G(\oo) \big(F(\oo+ k1_{[t,T]})- F(\oo)\big) \mu(\d\oo)\\
&\quad+\sum_{k\in \Z^d} \nu(k) \mathbb E\Big[G(X+k1_{[\tau, T]}) \big(F(X)- F(X+k1_{[\tau,T]})
\big)\Big]\\
&=: A_1+A_2. \end{split}\end{equation}By the definition of $\pi_k$, we have

\beg{equation*}\beg{split} A_2 &= \sum_{k\in \Z^d} \nu(k) \int_{M_T} G(\oo)
\mu_k(\d\oo)\int_0^T\big(F(\oo- k1_{[t,T]})-F(\oo)\big)\pi_k(\oo,\d t)\\
&= \int_{M_T} G(\oo) \bigg\{\sum_{k\in \Z^d} \ff{N_T^{(k)} (\oo)}T\bigg(\int_0^T F(\oo- k 1_{[t,T]})
\pi_k(\oo,\d t)-F(\oo)\bigg)\mu(\d\oo).\end{split}\end{equation*} 
Combining this with (\ref{2.32}), we arrive at 

\beq\label{G}-\E(F,G)=\int_{M_T} G(\oo) LF(\oo)\mu(\d\oo).\end{equation} This completes the proof. \end{proof} 

The formula (\ref{2.1}) indicates the transition rate of the associated jump process on $M_T$; that is, the process jumps from a state $\oo$
to $\oo +k1_{[t,T]}$ with rate $T^{-1} \nu(k)\d t$ while to $\oo- k_{[t,T]}$ with rate $T^{-1} N_T^{(k)}(\oo) \pi_k(\oo, \d t).$

In order to extend  these results for general $\nu$, we make use of an approximation procedure. To this end, let us first recall a standard construction of the Poisson process. Let $\scr N_t$ be the standard Poisson process on $\Z_+$ with intensity $\dd_1$ (i.e. the Dirac measure at $1$). Let $\tau_i$ be the $i$-th jump time of $\scr N_t$, i.e. 

$$\tau_1=\inf\{t\ge 0: \scr N_t>\scr N_{t-}\},\ \ \tau_{i}= \inf\{t\ge \tau_{i-1}:\ \scr N_t> 
\scr N_{t-}\}, i\ge 2.$$ Let $\{\xi_i\}_{i\ge 1}$ be i.i.d. sequence with distribution $\nu$ which are independent of $\scr N.$ 
Then

\beq\label{2.4} X_t:= \sum_{\tau_i\le t} \xi_i,\ \ \ t\ge 0\end{equation} is a Poisson process on $\R^d$ with intensity $\nu$. 

Now, for any $n\ge 1,$ let $\nu_n$ be the probability measure on $2^{-n}\Z^d$ with

$$\nu_n(2^{-n} k)=\nu(D_{n,k}),\ \ D_{n,k}:= \big\{x\in \R^d: x_i\in [2^{-n}k_i, 2^{-n}(k_i+1)), 1\le i\le d\big\},\ \ k\in \Z^d.$$ Let 

$$\xi_i^{(n)}= \sum_{k\in \Z^d} 2^{-n} k1_{D_{n,k}}(\xi_i),\ \ \ i\ge 1.$$  Let $\xi^{(n)}$ be determined by $\xi$ in the same way. Then 

$$X_t^{(n)}:= \sum_{\tau_i\le t} \xi_i^{(n)},\ \ \ t\ge 0$$ is a Poisson process with intensity $\nu_n.$ By Lemma \ref{L2.1}   we have 

$$\mathbb E F(X^{(n)}+ \xi^{(n)}1_{[\tau,T]})= \ff 1 T \mathbb E \big[F(X^{(n)})N_T\big]$$
for $F\in C_b(M_T).$ Letting $n\to\infty$ and using the dominated convergence theorem,
we prove (\ref{1.1}). The formula (\ref{G}) can be confirmed in the same way for 
$F,G\in C_b(M_T)$ and hence, also for $F,G\in \B_b(M_T)$ by the monotone class theorem. Therefore, we have proved the following result.

\beg{thm} \label{T2.2} $\mu$ is quasi-invariant under the map $\oo\mapsto \xi1_{[\tau,T]}$ such that $(\ref{1.1})$ holds. Consequently,  $(\E,\D(\E))$ is a conservative symmetric Dirichlet form on $L^2(M_T;\mu)$, and the generator  $(L,\D(L))$ with $\D(L)\supset \B_b(M_T)$ is given by $(\ref{2.1})$.\end{thm}

\section{The Poincar\'e and log-Sobolev inequalities}

We first   prove (\ref{P}). By the monotone class theorem and the fact that $\B_b(M_T)$ is dense in $\D(\E)$, it suffices to prove for the calss
$\scr FC_b$ of bounded cylindrical functions. 
\beq\label{2.3} \mu(F^2)\le \mu_2(\GG(F,F))+\mu(F)^2,\ \ \
F\in \scr FC_b.
\end{equation}
Let $\mathbb E^z $ be the expectation taking for the  the Poisson process $X$ starting at $z \in \R^d$. By taking $z =0$ (\ref{2.3}) follows from

\beq\label{2.4}\beg{split}  &\mathbb E^z  f^2(X_{t_1}, \cdots, X_{t_n}) -
\big( \mathbb E^z f (X_{t_1},\cdots, X_{t_n})\big)^2 \\
&\le  \mathbb E^z\int_{\R^d} \sum_{i=1}^n (t_i-t_{i-1})\big[
  f(X_{t_1}, \cdots, X_{t_n})\\
&\qquad\qquad\qquad-  f(X_{t_1}, \cdots,
X_{t_{i-1}}, X_{t_i} +x, \cdots, X_{t_n}+x)
\big]^2\nu(\d x),\ \Z  \in \R^d\end{split}\end{equation} for $0=t_0<t_1<\cdots<t_n\le T$ and $f\in C_0(\R^n).$
We shall prove this inequality by iterating in $n$.

(1) Let $n=1$ and $t_1=t\in (0,T].$ Then (\ref{2.4}) reduces to

\beq\label{2.5} \mathbb E^z f^2 (X_t)\le \big(\mathbb E^z f(X_t)\big)^2 +t \mathbb E^z \int_{\R^d} \big(f(X_t +x)-f(X_t)\big)^2\nu(\d x). \end{equation} Recall that the generator of
$X_t$ is (cf. \cite{O05})

 $$L_0 f( z)= \int_\R (f(z +x)- f(z ))\nu(\d x).$$ So  the associated square field is
 \beg{equation*}\beg{split} \GG_0(f,g)(z )&:= \big\{L_0(fg)- fL_0g -gL_0 f\big\}(z )\\
 & =\int_\R (f(z +x)- f( z))
 (g( z+x)-g(z ))\nu(\d x).\end{split}\end{equation*}
  Let $P_t^0 f(z )= \mathbb E^z f(X_t)$ be the corresponding Markov semigroup. We have

 \beq\label{2.6} P_t^0 f^2(z )- (P_t^0f( z))^2 =\int_0^t \ff{\d}{\d s} P_s^0(P_{t-s}^0 f)^2
 (z ) \d s = \int_0^t P_s^0 \GG_0(P_{t-s}^0f, P_{t-s}^0 f)(z )\d s.\end{equation}
 Since

 $$\mathbb E^z f(X_s+x)= \mathbb E^{z+x} f(X_s),\ \ \ s\ge 0, \ \  x\in \R^d,$$
 we have

 \beg{equation*}\beg{split} \GG_0(P_{t-s}^0 f, P_{t-s}^0 f))(z ) &= \int_{\R^d} \big(
 \mathbb E^z [ f(X_{t-s} +x)-f(X_{t-s})]\big)^2\nu(\d x)\\
 &\le \mathbb E^z \int_{\R^d} (f(X_{t-s} +x)-f(X_{t-s}))^2 \nu(\d x)= P_{t-s}
 \GG_0(f,f)(z).\end{split}\end{equation*}
 Then (\ref{2.6}) yields

 $$P_t^0 f^2(z ) - (P_t^0 f(z ))^2 \le \int_0^t P_s^0 P_{t-s}^0 \GG_0(f,f)(z )\d s= t P_t^0 \GG_0(f,f)( z).$$
 Thus, (\ref{2.5}) holds.

 (2) Assume that (\ref{2.4}) holds for $n=k,$ it remains to prove it for $n=k+1.$
 Let

 $$g(z ) = \mathbb E^z f( z , X_{t_2-t_1},\cdots, X_{t_{k+1}-t_1}),\ \ \Z  \in \R^d.$$
 By the assumption  we have

 \beg{equation*}\beg{split} &\mathbb E^z f^2( X_0, X_{t_2-t_1},\cdots, X_{t_{k+1}-t_1})
 -g(z )^2 \\
 &\le \mathbb E^z \int_{\R^d} \sum_{i=2}^{k+1} (t_i- t_{i-1}) \big[ f(X_0 ,
 X_{t_2-t_1},\cdots, X_{t_{k+1}-t_1})\\
 &\qquad\qquad- f(X_0 ,
 X_{t_2-t_1},\cdots, X_{t_{i-1}-t_1}, X_{t_i - t_1}+x, \cdots, X_{t_{k+1}-t_1}
 +x)\big]^2\nu(\d x),\ \ \ z\in \R^d.\end{split}\end{equation*}
Combining this with (\ref{2.5}) and using the Markov property,  we obtain

 \beg{equation*}\beg{split} &\mathbb E^z f^2(X_{t_1}, \cdots, X_{t_{k+1}}) =
 \mathbb E^z \mathbb E^{X_{t_1}} f^2(X_0, X_{t_2-t_1},\cdots, X_{t_{k+1}-t_1})\\
 &\le \mathbb E^z g^2(X_{t_1}) +\mathbb E^z  \mathbb E^{X_{t_1}}\int_{\R^d}
  \sum_{i=2}^{k+1} (t_i- t_{i-1}) \big[ f(X_0,
 X_{t_2-t_1},\cdots, X_{t_{k+1}-t_1})\\
 &\qquad\qquad- f(X_0,
 X_{t_2-t_1},\cdots, X_{t_{i-1}-t_1}, X_{t_i - t_1} +x, \cdots, X_{t_{k+1}-t_1}
 +x)\big]^2\nu(\d x )\\
 &\le (\mathbb E^z g(X_{t_1}))^2 + t_1 \mathbb E^z \int_{\R^d} (g(X_{t_1}+x)
 -g(X_{t_1}))^2\nu(\d  x)\\
 &\qquad + \mathbb E^z  \int_{\R^d} \sum_{i=2}^{k+1} (t_i- t_{i-1}) \big[ f(X_{t_1},
 \cdots, X_{t_{k+1}})\\
 &\qquad\qquad- f(X_{t_1},
 X_{t_2},\cdots, X_{t_{i-1}}, X_{t_i } +x, \cdots, X_{t_{k+1}}
 +x)\big]^2\nu(\d x )\\
 &= \big(\mathbb E^z f(X_{t_1},\cdots, X_{t_{k+1}})\big)^2 + \mathbb E^z  \int_{\R^d} \sum_{i=1}^{k+1} (t_i- t_{i-1}) \big[ f(X_{t_1},
 \cdots, X_{t_{k+1}})\\
 &\qquad\qquad- f(X_{t_1},
 X_{t_2},\cdots, X_{t_{i-1}}, X_{t_i } +x, \cdots, X_{t_{k+1}}
 +x)\big]^2\nu(\d x ).\end{split}\end{equation*} Therefore, (\ref{2.4}) holds
 for $n=k+1.$  

\

Finally, we intend to disprove the log-Sobolev inequality for any $C>0$. Let $\mu(F^2)=1.$ Noting that

$$\E(F,F) = \mathbb E (F(X+\xi 1_{[\tau,T]})-F(X))^2\le 2 \mathbb EF(X+\xi1_{[\tau,T]})^2 +2\mathbb E F(X)^2,$$ it follows from the definition of $\mu,\tt\mu$ and the formula 
(\ref{1.1}) that 

$$\E(F,F)\le 2\mu(F^2)+ \ff 2 T \mu(N_T F^2)\le 2 +\ll^{-1}  \mu(F^2\log F^2) +\ll^{-1}\log \mu(\e^{\ll N_T}).$$
Noting that under $\mu$ $N_T$ is a Poisson random variable  with intensity $T$, we conclude that $\mu(\e^{\ll N_T})<\infty$ for all $\ll>0$.  Thus, for any $\vv>0$ there exists
$C(\vv)>0$ such that 

$$\E(F,F)\le C(\vv) + \vv\mu(F^2\log F^2),\ \ \ \mu(F^2)=1.$$ Therefore, if (\ref{LS}) holds for some $C>0$, then there exists $C'>0$ such that 

$$\mu(F^2\log F^2)\le C',\ \ \ \mu(F^2)=1.$$ This is  wrong since the support of 
$\mu$ is not a finite set.

\end{document}